\documentclass[reqno]{amsart}
\usepackage{amsmath}
\usepackage{amssymb}
\usepackage{amsfonts}
\usepackage{euscript}
\usepackage[latin1]{inputenc}
\usepackage[T1]{fontenc}
\usepackage[pdftex]{hyperref}
\numberwithin{equation}{section}
\makeatletter\@addtoreset{equation}{section}

\DeclareMathSymbol{\subsetneqq}{\mathbin}{AMSb}{36}

\newcommand{\set}[1]{{\left\{{#1}\right\}}}

\newcommand{\scal}[1]{{\left\langle{#1}\right\rangle}}

\newcommand{\C}{\mathbb C}
\newcommand{\R}{\mathbb R}
\newcommand{\Z}{\mathbb Z}

\newcommand{\bz}{\bar z}
\newcommand{\X}{\xi}

\newcommand{\Nu}{\sigma}

\newcommand{\La}{\mathfrak{L}^{\nu,\mu}_\tau}
\newcommand{\Laa}{\mathfrak{L}^{\Nu}}

\newtheorem {theorem}{Theorem}[section]

\newtheorem {lemma}[theorem]{Lemma}
\newtheorem {proposition}[theorem]{Proposition}

\newtheorem {corollary}[theorem]{Corollary}

\newtheorem {claim}[theorem]{Claim}

\begin{document}

\title[Magnetic Schr\"odinger operator on mixed automorphic forms]
{On spectral analysis of a magnetic Schr\"odinger operator on planar mixed automorphic forms}
\author{A. Ghanmi}
\address{Department of Mathematics, Faculty of  Sciences, P.O. Box 1014\newline
   Mohammed V University,  Agdal,  10 000 Rabat - Morocco  }
   \email{allalghanmi@gmail.com}

\date{\today \newline {\it 2000 Mathematics Subject Classification:} 58C40; 11F03; 32N05; 30C40; 81Q10 \newline
}

\maketitle

\begin{abstract} \quad
We characterize the space of the so-called planar mixed automorphic forms of
type $(\nu,\mu)$ with respect to an equivariant pair $(\rho,\tau)$ as the image of the usual
 automorphic forms by an appropriate transform and we investigate some concrete basic spectral properties of a magnetic
 Schr\"odinger operator acting on them.
The associated polynomials constitute  classes of generalized complex polynomials of
Hermite type.
\end{abstract}

 \section{Introduction}
Let $H_{\theta,V}$ be the Schr\"odinger operator
\begin{eqnarray}\label{SchrOp}
H_{\theta,V} f := - \Delta f - 2i (df|\theta) + \big(i d^{*}\theta +  |\theta|^2 +V \big) f
 \end{eqnarray}
on a Riemannian manifold $M$, where $\Delta$ is the Laplace-Beltarmi operator, $i=\sqrt{-1} $, $\theta$ is a real $\mathcal{C}^\infty$ differential $1$-form, $V$ is a real valued
$\mathcal{C}^\infty$ function and $(\cdot|\cdot)$ is the Hermitian inner product induced from the Riemannian metric.
Such operators, for both compact or noncompact manifold, arise in many problems both of classical
mechanics and mathematics and modern mathematical physics; see \cite{Shigekawa87,CFKS87,Agmon95,Shubin01} and references therein. For $M$ being the Euclidean plane $\R^2=\C$, $\theta$ being
  the canonical vector potential $\theta=i\nu(\bar z dz - zd\bar z) $; $\nu\in\R$, and $V=0$,
  the operator \eqref{SchrOp} leads to the special Hermite operator $\mathfrak{L}^\nu$ \cite{Thangavelu93,Wong98},
\begin{equation}\label{LandauHamiltonian}  \mathfrak{L}^\nu =   -
\frac{\partial^2}{\partial z \partial \bar z} - \nu
(z\frac{\partial}{\partial z}-\bar z \frac{\partial}{\partial \bar
z}) + \nu^2 |z|^2 ,\end{equation}
which  describes in physics the quantum behavior of a charged particle confined in the
plane, under the action of a constant magnetic field.
The spectral theory of  $\mathfrak{L}^\nu $ when acting on the Hilbert space $L^2(\C; d\lambda)$, of square integrable functions on $\C$ w.r.t. the usual Lebesgue measure $d\lambda$, is well known (see for example \cite{CFKS87,Matsumoto96,AIM00,FerapontovVeselov01,GI05}). For $\mathfrak{L}^\nu$ on the space of Landau automorphic functions of magnitude $\nu$, a systematic study is recently presented in \cite{GI-JMP08}.
In this paper, we shall consider a special Schr\"odinger operator acting on the so-called mixed automorphic forms.
Thus, let $\mathbf{T}=\set{\lambda\in \C; |\lambda|=1}$
and consider the semidirect product group
\[G=\mathbf{T}\rtimes \C =\set{g=\left(\begin {array} {c c} a & b \\ 0 & 1
\end {array} \right )=:[a,b]; \quad a\in \mathbf{T} , ~ b\in \C}\]
 operating on $\C$ by
 $g\cdot z = a z +b$.
  Define $j^\alpha$; $\alpha\in \R$, to be the automorphic factor
given by
\begin{equation}\label{AutomFactorExplicit}
 j^\alpha(g,z)=e^{2i\alpha  \Im \scal{z,g^{-1}\cdot 0}}, \quad g \in G, z\in \C,
\end{equation}
where here and elsewhere $\Im z$ denotes the imaginary part of the complex number $z$ and
$\scal{\cdot,\cdot}$ denotes the usual hermitian scalar product on $\C$.
Also, let  $\Gamma$ be a given uniform lattice in $\C$ and  $(\rho,\tau)$ be an equivariant pair
\cite{SatakeBook80,HuntMeyer85,Abdulali95,LeeLectNotes04}.
That is,  $\rho$ is a $G$-endomorphism  and $\tau : \C \to \C$ a smooth compatible mapping  such that
\begin{equation}
\label{EquivCond} \tau(g.z)= \rho(g)\cdot\tau(z)
\end{equation}
 for every $g\in G$ and $z\in\C$.
Associated to given real numbers $\nu, \mu >0$ we perform
the space $\mathcal{M}^{\nu,\mu}_{\tau}(\C)$ of mixed $\Gamma$-automorphic forms of type $(\nu,\mu)$ w.r.t $(\rho,\tau)$, i.e., the space of all  $\mathcal{C}^\infty$ complex-valued functions $F$ on $\C$
satisfying the functional equation
 \begin{equation}\label{FunctionalEq1}
 F(\gamma\cdot z) =  j^{-\nu}(\gamma, z)j^{-\mu}(\rho(\gamma), \tau(z)) F(z)
\end{equation}
for every $\gamma\in \Gamma$ and $z\in \C$. Such notion of mixed automorphic forms has been introduced by Hunt and Meyer
\cite{HuntMeyer85} and extensively studied by M.H. Lee (see \cite{LeeLectNotes04} and the references therein).
It generalizes in fact the one of automorphic forms (take for example $\mu=0$ or let $\rho$ and $\tau$
be the identity maps) and appears essentially in the context of number theory and
algebraic geometry.

\noindent The main aim of this paper is firstly to characterize the space $\mathcal{M}^{\nu,\mu}_{\tau}(\C)$ and show that
the two pictures of mixed autmorphic forms and classical automorphic forms can be connected by a
special transform, and secondly to investigate some concrete spectral properties of the appropriate magnetic Schr\"odinger operator
$H_{\theta_\tau^{\nu,\mu},0}$, corresponding to
\begin{equation}\label{ThetaForm1}
\theta= \theta^{\nu,\mu}_\tau(z) := i \Big\{\nu(\bar z dz -  zd\bar z) +\mu(\overline{\tau}d\tau - \tau d\overline{\tau}) \Big\},
\end{equation}  when acting on the free Hilbert space $L^2(\C; d\lambda)$ and the space $\mathcal{M}^{\nu,\mu}_{\tau} (\C)$ of mixed automorphic forms \eqref{FunctionalEq1}.

This paper is organized as follows. In Section 2,
we investigate some standard basic properties of the invariant magnetic Schr\"odinger operator $\La$. In Section 3, we give a characterization of
$\mathcal{M}^{\nu,\mu}_{\tau}(\C)$ in terms the usual automorphic functions (Theorem \ref{Thm2}).
Section 4 is devoted to describe the structure of the associated eigenspaces
 in $L^2(\C;d\lambda)$ and $\mathcal{M}^{\nu,\mu}_{\tau}(\C)$ for general equivariant map $\tau$.
 The particular case of $\tau=\tau_{_{\mathfrak{h}}}$ being the one associated the inner $G$-automorphism $\rho_{_{\mathfrak{h}}}(g)
      :=\mathfrak{h}g \mathfrak{h}^{-1}$; $\mathfrak{h}\in G$, is considered for illustration.
      Section 5 deals with a discussion of the derived classes of complex polynomials of Hermite type.

\section{The magnetic Schr\"odinger operator $\La$}

Keep notation as in the previous section and set $\La  :=(1/4) H_{\theta^{\nu,\mu}_\tau,0}$,
\begin{eqnarray}\label{MagnSchrOp}
 \La f := \frac 14 \bigg\{-\Delta f - 2i (df|\theta^{\nu,\mu}_\tau)
 + \bigg(i d^{*}\theta^{\nu,\mu}_\tau +  |\theta^{\nu,\mu}_\tau|^2 \bigg) f\bigg\},
 \end{eqnarray}
 where $\Delta$ is the Euclidean Laplacian on $\C$, $\theta^{\nu,\mu}_\tau$ is the differential $1$-form, \eqref{ThetaForm1},
\begin{equation}\label{ThetaForm}
\theta^{\nu,\mu}_\tau(z) : =i(\overline{S^{\nu,\mu}_\tau}dz - S^{\nu,\mu}_\tau d\bar z)
\end{equation}
with $S^{\nu,\mu}_\tau$ stands for the complex-valued function
\begin{equation}\label{ThetaForm3}
S^{\nu,\mu}_\tau (z)= \nu z +\mu \Big(\tau \frac{\partial\bar
\tau}{\partial  \bar z}- \bar \tau \frac{\partial \tau}{\partial
\bar z} \Big).
\end{equation}

 The main observation of this section is included in the following result.

\begin{proposition}\label{PropMagnSchrOp}  The operator $\La$ is a Schr\"odinger operator with constant magnetic field whose the
explicit expression in terms of ordinary partial derivatives is given by
\begin{equation}
\La   = - \frac{\partial^2 }{\partial z\partial \bar z}    - \Big(S^{\nu,\mu}_\tau
\frac{\partial }{\partial z} -
\overline{S^{\nu,\mu}_\tau} \frac{\partial }{\partial \bar z} \Big)
+ | S^{\nu,\mu}_\tau |^2   - \frac \mu{4}(\tau\Delta \bar \tau -\bar \tau \Delta \tau). \label{ExplicitAutoLap}
\end{equation}
Moreover, $\La $ satisfies the following supersymmetric relations
   \begin{eqnarray}
                   \widetilde{A^{\nu,\mu}_\tau} A^{\nu,\mu}_\tau  + B^{\nu,\mu}_\tau = \La , \quad \mbox{and} \quad
                     A^{\nu,\mu}_\tau \widetilde{A^{\nu,\mu}_\tau} - B^{\nu,\mu}_\tau = \La , \label{Super1}
   \end{eqnarray}
where $A^{\nu,\mu}_\tau$ and $\widetilde{A^{\nu,\mu}_\tau}$ are the
 first order differential operators
\begin{equation}\label{AnniCrea}A^{\nu,\mu}_\tau= \frac{\partial}{\partial \bar z} + S^{\nu,\mu}_\tau \quad \mbox{and} \quad
\widetilde{A^{\nu,\mu}_\tau}=  -\frac{\partial}{\partial z} + \overline{S^{\nu,\mu}_\tau} ,\end{equation}
and $B^{\nu,\mu}_\tau $ is the real valued constant given by
$B^{\nu,\mu}_\tau= \nu +\mu\big(\big|\frac{\partial \tau}{\partial z}\big|^2- \big|\frac{\partial \tau}{\partial\bar z}\big|^2\big).$
\end{proposition}

\begin{proof}
For general vector potential $\theta$, one can show that the derived magnetic field
 $\mathcal{B} := d\theta$ is uniform if and only if $\theta$ and its
 pullback $g^{*}\theta$ by the holomorphic mapping
 $z \mapsto g\cdot z$ belong to the same de Rham co-homology group. This is clearly satisfied when $\theta=\theta^{\nu,\mu}_\tau$,
but here we will provide a direct and simple proof. In fact,
the derived magnetic field is given by $\mathcal{B}^{\nu,\mu}_\tau(z)=-2iB^{\nu,\mu}_\tau(z)dz\wedge d\bar z $, where  \begin{equation}\label{MagnField}
B^{\nu,\mu}_\tau(z)= \nu +\mu\bigg(\bigg|\frac{\partial \tau}{\partial z}({z})\bigg|^2- \bigg|\frac{\partial \tau}{\partial\bar z}({z})\bigg|^2\bigg).
\end{equation}
Next by writing the $G$-endomorphism $\rho$ as $\rho(g)=[\chi(g),\psi(g)]\in G=\mathbf{T}\rtimes \C$
and differentiating the equivariant condition $ \tau(g\cdot z) =\rho(g)\cdot \tau(z)= \chi(g)\tau(z) + \psi(g),$
 we get
$$ \frac{\partial \tau}{\partial z}({g\cdot z})= \chi(g) \overline{\big(\frac{g\cdot z}{\partial z}\big)}\frac{\partial \tau}{\partial z}(z) \quad \mbox{and} \quad
   \frac{\partial \tau}{\partial \bar z}({g\cdot z})= \chi(g)\big(\frac{g\cdot z}{\partial z}\big)\frac{\partial \tau}{\partial \bar z}(z).$$
   Hence, since $(\frac{g\cdot z}{\partial z})\chi(g) $ belongs to $\mathbf{T}$,
 we see that the function $z \mapsto B^{\nu,\mu}_\tau(z)$ is $G$-invariant 
  and therefore is constant on $\C$.

Now to get the explicit expression of $\La$, let recall first that the adjoint operation is taken w.r.t. the Hermitian
scalar product on compactly supported differential forms
$ (\alpha,\beta):=\int_{\C}\alpha\wedge\star\beta, $
 where $\star $ is the Hodge star operator canonically associated with the Euclidean metric of $\C$. Thus, the adjoint $d^*$ of $d$
and $\theta^{*}$ of $\theta$ are given respectively by $d^*=-\star d\star$
and $\theta^{*}= \star \theta \wedge \star .$ Therefore, one can check by direct computation the following
\begin{align*}
 & i) \qquad    (df|\theta^{\nu,\mu}_\tau) = -2i \bigg(S^{\nu,\mu}_\tau \frac{\partial f}{\partial  z}  -
                     \overline{S^{\nu,\mu}_\tau} \frac{\partial f}{\partial
                    \bar z}\bigg) ,\\ 
 & ii) \qquad   d^{*}(\theta^{\nu,\mu}_\tau) = 2i \bigg( \frac{\partial S^{\nu,\mu}_\tau}{\partial  z}  -
                      \frac{\partial \overline{S^{\nu,\mu}_\tau}}{\partial
                     \bar z}\bigg) ,\\  
 & iii) \qquad  |\theta^{\nu,\mu}_\tau |^2 = 4| S^{\nu,\mu}_\tau |^2,
\end{align*}
Note finally that \eqref{Super1} can be handled  by easy direct computation.
\end{proof}


As immediate consequence of Proposition \ref{PropMagnSchrOp}, we get the following

 \begin{lemma}\label{Gauge}  There exists a real-valued function $\varphi^{\nu,\mu}_\tau$, with $\varphi^{\nu,\mu}_\tau(0)=0$, such that
  the annihilator and creator operators $A^{\nu,\mu}_\tau$ 
 and $\widetilde{A^{\nu,\mu}_\tau}$ can be rewritten as
 \begin{align}   A^{\nu,\mu}_\tau  & = e^{-B |z|^2 -i \varphi^{\nu,\mu}_\tau(z)} \frac{\partial}{\partial \bar z} e^{B |z|^2 +i \varphi^{\nu,\mu}_\tau(z)} .\label{Super4}\\
 \widetilde{A^{\nu,\mu}_\tau}  & = e^{B |z|^2 -i \varphi^{\nu,\mu}_\tau(z)} \frac{\partial}{\partial z} e^{-B |z|^2 +i \varphi^{\nu,\mu}_\tau(z)} .\label{Super5}
   \end{align}
 \end{lemma}

\begin{proof} Since $B^{\nu,\mu}_\tau$ is a constant, one checks that $d\theta^{\nu,\mu}_\tau = d\theta^{B^{\nu,\mu}_\tau},$ where $\theta^{B^{\nu,\mu}_\tau}:=iB^{\nu,\mu}_\tau (\bar z dz-zd\bar z)$.
 Therefore, there exists a function
$\widetilde{\varphi^{\nu,\mu}_\tau}: \C \to \C$ such that $\theta^{\nu,\mu}_\tau = \theta^{B^{\nu,\mu}_\tau}+ d \widetilde{\varphi^{\nu,\mu}_\tau}.$
Hence, $\widetilde{\varphi^{\nu,\mu}_\tau}$ is solution of the following system of first order partial differential equations
\begin{eqnarray}\label{Existence1stDiffEq}
\left\{ \begin{array}{lll}
\frac{\partial \widetilde{\varphi^{\nu,\mu}_\tau}}{\partial z} = i\big(\overline{S^{\nu,\mu}_\tau }- B^{\nu,\mu}_\tau \bar z\big)\\
\quad \\
\frac{\partial \widetilde{\varphi^{\nu,\mu}_\tau}}{\partial \bar z} = -i\big(S^{\nu,\mu}_\tau - B^{\nu,\mu}_\tau z\big)
\end{array} \right.
\end{eqnarray}
or equivalently
\begin{eqnarray}\label{Existence1stDiffEqEq}
\left\{ \begin{array}{lll}
\frac{\partial \widetilde{\varphi^{\nu,\mu}_\tau}}{\partial \bar z} = -i\big(S^{\nu,\mu}_\tau - B^{\nu,\mu}_\tau z\big)\\
\quad \\
\Im\widetilde{\varphi^{\nu,\mu}_\tau}= {Constant}
\end{array} \right. .
\end{eqnarray}
This assures the existence of solutions of the following first order differential equation
\begin{equation}  \frac{\partial \phi}{\partial \bar z} =   S^{\nu,\mu}_\tau
= \nu z +\mu( \tau \frac{\partial \bar\tau}{\partial \bar z}- \bar \tau\frac{\partial \tau}{\partial\bar z})
\label{1stDiffEq}
   \end{equation}
   and therefore of the one arising from
   \begin{equation*}  A^{\nu,\mu}_\tau
= e^{-\phi} \frac{\partial }{\partial \bar z}e^{\phi}
   \end{equation*}
Indeed, solutions of \eqref{1stDiffEq} are of the form \begin{equation}
\phi(z)=  B|z|^2 + i\widetilde{\varphi^{\nu,\mu}_\tau}(z) + h(z), \label{SolFormDiffEq}
   \end{equation}
    where $h$ is any arbitrary holomorphic function. This finishes the proof by choosing $h\equiv 0$ and taking $\varphi^{\nu,\mu}_\tau:=\Re (\widetilde{\varphi^{\nu,\mu}_\tau}-\widetilde{\varphi^{\nu,\mu}_\tau}(0))$.
    \end{proof}

    We conclude this section by giving an invariance property for $\La$ by
 the unitary transformations $\mathcal{T}^{\nu,\mu}_g$; $g\in G$, defined by
\begin{equation}\label{Transformation}
 [\mathcal{T}^{\nu,\mu}_gf](z):= J^{\nu,\mu}_{\rho,\tau}
(g,z) f(g\cdot z) ,
\end{equation}
where $J^{\nu,\mu}_{\rho,\tau}$ is the complex valued mapping on  $ G\times \C$ given by
\begin{equation}\label{AutomFactor}
J^{\nu,\mu}_{\rho,\tau} (g,z):= j^\nu(g, z) j^\mu(\rho(g), \tau(z)).
\end{equation}
 Namely, we assert the following

\begin{proposition}\label{PropMagnSchrOp2}
 The magnetic Schr\"odinger operator $\La$ can be realized in terms of the deformed de Rham differential
 $\nabla^{\nu,\mu}_\tau f := df+if\theta^{\nu,\mu}_\tau$ as
$\La  = (\frac 14)(\nabla^{\nu,\mu}_\tau )^{*}\nabla^{\nu,\mu}_\tau .$
 Moreover, it satisfies the invariance property
$ \mathcal{T}^{\nu,\mu}_g\La  = \La \mathcal{T}^{\nu,\mu}_g$ for all $g\in G$.
\end{proposition}

\begin{proof}
We have
\begin{align*}
\nabla^{*}\nabla f &= d^{*}d f  + i[d^{*}(f \theta) - \theta^{*}(df) ]  + \theta^{*}(f\theta) \\
                   &= -\Delta  - i[\star ((d \bar f) \wedge (\star\theta))   + \theta^{*}(df) ]  + \big(i d^{*}(\theta) + |\theta|^2\big)f \\
                   &= -\Delta  - i\star [(d \bar f) \wedge \star\theta   + \theta\wedge \star df ]  + \big(i d^{*}(\theta) + |\theta|^2\big)f \\
                   &= -\Delta  - 2i (df|\theta)   + \big(i d^{*}(\theta) + |\theta|^2\big)f =  4\La .
\end{align*}
The second equation follows using the facts $d^{*}d = \Delta$, $d^{*}(f \theta) = -\star (d \bar f \wedge (\star\theta))  + f d^{*}(\theta)$
and $\theta^{*}(f\theta)=|\theta|^2f$, while the third one is due to $\star[(d \bar f) \wedge (\star\theta)   + \theta\wedge \star df]
= 2\star ((d \bar f) \wedge (\star\theta) ) = (df) \wedge (\star\theta) 
$ which follows from  $\alpha\wedge \star \beta = \bar \beta \wedge \star \bar\alpha$, keeping in mind that $\theta$ is real, combined with $\star (\alpha\wedge \star \beta) = (\beta|\alpha)$.

For the invariance property of $\La$ by the transformations $\mathcal{T}^{\nu,\mu}_g$, let mention first that the differential $1$-form $\theta^{\nu,\mu}_\tau$
satisfies the following
\begin{equation}\label{CoHom}
 [g^{*}\theta^{\nu,\mu}_\tau](z)
= \theta^{\nu,\mu}_\tau(z)  -i \frac{d(J^{\nu,\mu}_{\rho,\tau} (g,z))}{J^{\nu,\mu}_{\rho,\tau} (g,z)}.
\end{equation}
This holds by direct computation and the use of the equivariant condition \eqref{EquivCond}. Now,
 using the well known facts $g^{*}d=dg^{*}$ as well as $g^{*}(\alpha\wedge
\beta)= g^{*}\alpha\wedge g^{*}\beta,$ together with the established fact \eqref{CoHom},
 one obtains
\begin{align*}
\mathcal{T}^{\nu,\mu}_g\big(\nabla^{\nu,\mu}_\tau f\big)
 &= J^{\nu,\mu}_{\rho,\tau} (g,z) \Big[g^{*}\Big((d+i\theta^{\nu,\mu}_\tau)f\Big)\Big]\\
&= J^{\nu,\mu}_{\rho,\tau} (g,z) \Big(d[g^{*}f]+i[g^{*}\theta^{\nu,\mu}_\tau ]\wedge [g^{*}f]\Big)\\
&= J^{\nu,\mu}_{\rho,\tau} (g,z) d [g^{*}f]+ d(J^{\nu,\mu}_{\rho,\tau} (g,z)) [g^{*}f] +
i J^{\nu,\mu}_{\rho,\tau} (g,z)\theta^{\nu,\mu}_\tau [g^{*}f]\\
&= d(J^{\nu,\mu}_{\rho,\tau} (g,z) [g^{*}f])
 + i\theta^{\nu,\mu}_\tau J^{\nu,\mu}_{\rho,\tau} (g,z)[g^{*}f]\\
& = (d+i\theta^{\nu,\mu}_\tau )
 \mathcal{T}^{\nu,\mu}_g f = \nabla^{\nu,\mu}_\tau \big(\mathcal{T}^{\nu,\mu}_g f\big).
\end{align*}
Hence, since $\mathcal{T}^{\nu,\mu}_g$ commutes also with
$(\nabla^{\nu,\mu}_\tau)^{*}$ for $\mathcal{T}^{\nu,\mu}_g$
being an unitary transformation, we deduce easily from
$\La=(1/4)\nabla^{*}\nabla$ that $\La$ and $\mathcal{T}^{\nu,\mu}_g$ commute. This ends the proof.
\end{proof}

\section{A characterization of the space $\mathcal{M}^{\nu,\mu}_{\tau}(\C)$}

 Let $\phi^{\nu,\mu}_\rho(g,g')$ be the real-valued function defined on $G\times G$ by
\begin{equation}\label{PhaseFactor}
\phi^{\nu,\mu}_\rho(g,g'):= \Im  \Big(\nu \scal{g^{-1}\cdot
0,g'\cdot 0} + \mu \scal{\rho({g^{-1}})\cdot
0,\rho(g')\cdot 0} \Big)
\end{equation}
and recall that $J^{\nu,\mu}_{\rho,\tau}$, in \eqref{AutomFactor}, is given by
$J^{\nu,\mu}_{\rho,\tau} (g,z):= j^\nu(g, z) j^\mu(\rho(g), \tau(z)). $
Then, we have

\begin{lemma}
The mapping $J^{\nu,\mu}_{\rho,\tau}$ satisfies the chain rule
\begin{align}
J^{\nu,\mu}_{\rho,\tau} (gg',z) = e^{2i\phi^{\nu,\mu}_\rho(g,g')}
 J^{\nu,\mu}_{\rho,\tau} (g,g'\cdot z)
J^{\nu,\mu}_{\rho,\tau} (g',z).  \label{ChainRule}
\end{align}
\end{lemma}

\begin{proof} For every $g,g'\in G$ and $z\in \C$, we have
\begin{align*}
J^{\nu,\mu}_{\rho,\tau} (gg',z) = j^{\nu} (gg',z)j^{\mu}(\rho(gg'),\tau(z)) = j^{\nu} (gg',z)j^{\mu}(\rho(g)\rho(g'),\tau(z)).
\end{align*}
Next, by applying the fact that the automorphic factor $j^{\alpha}(\cdot,\cdot)$ satisfies
$$ j^{\alpha}(hh',w) =  e^{2i\alpha \Im  \scal{h^{-1}\cdot
0,h'\cdot 0}} j^{\alpha}(h, h'\cdot w) j^{\alpha}(h',w)$$ successively for $h=g, h'=g'$ and $h=\rho(g), h'=\rho(g')$, we deduce
\begin{align*}
J^{\nu,\mu}_{\rho,\tau} (gg',z)
&= e^{2i\phi^{\nu,\mu}_\rho(g,g')} j^{\nu}(g, g'\cdot z)  j^{\mu}(\rho(g), \rho(g')\cdot \tau(z)) J^{\nu,\mu}_{\rho,\tau} (g',z)\\
&= e^{2i\phi^{\nu,\mu}_\rho(g,g')} j^{\nu}(g, g'\cdot z)  j^{\mu}(\rho(g), \tau(g'\cdot z)) J^{\nu,\mu}_{\rho,\tau} (g',z).
\end{align*}
The last equation follows making use of the equivariant condition \eqref{EquivCond}. 
\end{proof}

\noindent According to the above result, the unitary transformations $\mathcal{T}^{\nu,\mu}_g$; $g\in G$, \eqref{Transformation},
 define a projective representation of the group $G$ on the space of $\mathcal{C}^\infty$
functions on $\C$. Therefore,
we can check the following statement on the nontriviality of the vector space of mixed $\Gamma$-automorphic forms of type $(\nu,\mu)$,
\begin{equation}\label{MixedForms}
\mathcal{M}^{\nu,\mu}_{\tau}(\C):= \set{F: \C \stackrel{\mathcal{C}^\infty}{\longrightarrow} \C;  \quad  F(z+\gamma) =
J^{-\nu,-\mu}_{\rho,\tau}(\gamma,z) F(z) }.
\end{equation}

\begin{lemma}\label{NonTrivial}
The functional space $\mathcal{M}^{\nu,\mu}_{\tau}(\C)$ is nontrivial if and only if the real-valued function $(1/\pi)\phi^{\nu,\mu}_\rho$ in \eqref{PhaseFactor} takes integer values on $\Gamma\times \Gamma$. In this case, $\mathcal{M}^{\nu,\mu}_{\tau}(\C)$
can be realized as the space of cross sections on a line bundle over
the complex torus $\C/\Gamma$.
\end{lemma}

\begin{proof} Making use of \eqref{ChainRule} and the equivariant condition (\ref{EquivCond}), the proof can be handled along the lines of the proof of \cite[Proposition 3.1]{GI-JMP08} (see also Remark 3.3 there).
\end{proof}

Now in order to state the main result of this paper, Theorem \ref{Thm2} below, let consider the transformation
\begin{eqnarray}\label{GaugeTransformation}
[\mathcal{W}^{\nu,\mu}_\tau  (f)](z):= e^{i \varphi^{\nu,\mu}_\tau(z)}f(z),
\end{eqnarray}
where $\varphi^{\nu,\mu}_\tau$ is the real-valued function considered in Lemma \ref{Gauge}. Then, we have the following

 \begin{theorem} \label{Thm2}
  The image of $\mathcal{M}^{\nu,\mu}_{\tau}(\C)$ by the special transformation \eqref{GaugeTransformation}
  is the space of Landau $(\Gamma,\chi_\tau)$-automorphic functions of magnitude $B^{\nu,\mu}_\tau= \nu +\mu(|\frac{\partial \tau}{\partial z}|^2- |\frac{\partial \tau}{\partial\bar z}|^2)$ and multiplier
 $$\chi_\tau(\gamma)=\exp\big({2i\varphi^{\nu,\mu}_\tau(\gamma)-2i\mu  \Im \scal{\tau(0),\rho(\gamma)^{-1}\cdot 0}}\big).$$
  That is
  $$\mathcal{W}^{\nu,\mu}_\tau (\mathcal{M}^{\nu,\mu}_{\tau}(\C) )
  =\set{F:\C \stackrel{\mathcal{C}^\infty}{\longrightarrow} \C;  \quad F(z+\gamma)= \chi_\tau(\gamma)j^{-B^{\nu,\mu}_{\tau}}(\gamma, z)  F(z)}
 .$$
\end{theorem}

   \begin{proof}
  We have to prove that $\mathcal{W}^{\nu,\mu}_\tau F$ belongs to
  $$\mathcal{F}^{B^{\nu,\mu}_{\tau}}_{\Gamma,\chi_{\tau}}(\C) : =\set{F ;\, \mathcal{C}^\infty, ~ F(z+\gamma)=
\chi_\tau(\gamma)j^{-B^{\nu,\mu}_{\tau}}(\gamma, z)  F(z)}$$
   whenever  $F \in \mathcal{M}^{\nu,\mu}_{\tau}(\C)$, where $\chi_\tau(\gamma)=\exp\big(2i\varphi^{\nu,\mu}_\tau(\gamma)-2i\mu  \Im \scal{\tau(0),\rho(\gamma)^{-1}\cdot 0}\big).$
Indeed, we have
\begin{align*}
[\mathcal{W}^{\nu,\mu}_\tau F](z+\gamma)
& := e^{i\varphi^{\nu,\mu}_\tau(z+\gamma)} F(z+\gamma) \\
&  = e^{i\varphi^{\nu,\mu}_\tau(z+\gamma)} j^{-\nu}(\gamma, z)j^{-\mu}(\rho(\gamma), \tau(z)) F(z) \\
&  = e^{i(\varphi^{\nu,\mu}_\tau(z+\gamma)-\varphi^{\nu,\mu}_\tau(z))} j^{-\nu}(\gamma, z)j^{-\mu}(\rho(\gamma), \tau(z)) [\mathcal{W}^{\nu,\mu}_\tau F](z) \\
&  = \widehat{\chi_\tau}(z;\gamma)  j^{-B^{\nu,\mu}_\tau}(\gamma, z) [\mathcal{W}^{\nu,\mu}_\tau F](z) ,
\end{align*}
where we have set
\begin{equation}\label{MyPseudoChi}
\widehat{\chi_\tau}(z;\gamma):=  e^{i(\varphi^{\nu,\mu}_\tau(z+\gamma)- \varphi^{\nu,\mu}_\tau(z))} e^{-2i([B^{\nu,\mu}_{\tau}-\nu] \Im\scal{z,\gamma} +\mu  \Im \scal{\tau(z),\rho(\gamma)^{-1}\cdot 0})}.
\end{equation}
Thus, we claim

 \begin{claim}\label{LemTransMixed}
 The function $\widehat{\chi_\tau}$ as defined by \eqref{MyPseudoChi} is independent of the variable $z$.
 \end{claim}

\noindent Therefore, we get
 $$ [\mathcal{W}^{\nu,\mu}_\tau F](z+\gamma) = \chi_\tau(\gamma)  j^{-B^{\nu,\mu}_\tau}(\gamma, z) [\mathcal{W}^{\nu,\mu}_\tau F](z)$$
 with $\chi_\tau(\gamma):=\widehat{\chi_\tau}(0;\gamma) .$ The proof is completed
 \end{proof}

 \begin{proof}[Proof of Claim \ref{LemTransMixed}]
 Differentiation of $\widehat{\chi_\tau}(z;\gamma)$ w.r.t. variable $z$  gives
\begin{align} \frac{\partial \widehat{\chi_\tau}}{\partial z}   & =  i \bigg(
 \frac{\partial \varphi^{\nu,\mu}_\tau}{\partial z}(z+\gamma) - \frac{\partial\varphi^{\nu,\mu}_\tau}{\partial z}(z)\bigg) \widehat{\chi_\tau} \label{Indep1} \\ & \qquad
 - \bigg([B^{\nu,\mu}_{\tau}-\nu] \bar\gamma  + \mu  \bigg[ \overline{a_\gamma} \frac{\partial\tau}{\partial z}(z) - a_\gamma \frac{\partial\bar \tau}{\partial z}(z)\bigg]
 \bigg)\widehat{\chi_\tau}, \nonumber\end{align}
 where we have set $a_\gamma=\rho(\gamma^{-1})\cdot 0$. In the other hand,
 since $$ \frac{\partial \varphi^{\nu,\mu}_\tau}{\partial z} =\frac{\partial \widetilde{\varphi^{\nu,\mu}_\tau}}{\partial z}  ,$$
 one deduces from \eqref{Existence1stDiffEq} together  with
 the equivariant condition $\tau(z+\gamma)= \rho(\gamma)\cdot \tau(z)$, the following
 \begin{align} i\bigg(\frac{\partial \varphi^{\nu,\mu}_\tau}{\partial z}(z+\gamma)- \frac{\partial\varphi^{\nu,\mu}_\tau}{\partial z}(z)\bigg)
& =   B^{\nu,\mu}_{\tau}\bar\gamma +  \overline{S^{\nu,\mu}_{\tau}}(z) - \overline{S^{\nu,\mu}_{\tau}}(z+\gamma) \nonumber\\
& =  [B^{\nu,\mu}_{\tau}-\nu]\bar\gamma   + \mu\bigg[ \overline{a_\gamma} \frac{\partial\tau}{\partial z}(z) - a_\gamma \frac{\partial\bar \tau}{\partial z}(z)\bigg] .   \label{Indep2}
   \end{align}
  Thus from \eqref{Indep1} and \eqref{Indep2}, we conclude that $\frac{\partial \widehat{\chi_\tau}}{\partial z}=0$. Similarly, one obtains also $\frac{\partial \widehat{\chi_\tau}}{\partial \bar z}=0$. This ends the proof of Claim \ref{LemTransMixed}.
 \end{proof}

 According to the results of Lemma \ref{NonTrivial} and Theorem \ref{Thm2} above together with Proposition 2.1 of \cite{GI-JMP08}, we deduce easily the following

 \begin{corollary} \label{CorRDQ}
 The function $\chi_\tau(\gamma)=\exp\big(2i\varphi^{\nu,\mu}_\tau(\gamma)-2i\mu  \Im \scal{\tau(0),\rho(\gamma)^{-1}\cdot 0}\big)$
 satisfies the pseudo-character property
 $$\chi_\tau(\gamma+\gamma')=e^{2i B^{\nu,\mu}_\tau \Im \scal{\gamma,\gamma'} } \chi_\tau(\gamma)\chi_\tau(\gamma')$$
 if and only if  $(1/\pi)\phi^{\nu,\mu}_\rho$ in \eqref{PhaseFactor} takes integer values on $\Gamma\times \Gamma$.
 \end{corollary}

\section{ On the concrete spectral theory of $\La$ on $L^2(\C;d\lambda)$ and ${\mathcal{M}}^{\nu,\mu}_\tau(\C)$}

  Note first the the spectrum of $\La$ is purely discrete and consists of infinitely degenerate eigenvalues (Landau levels)
 $$E_k= B(2k+1)=(2k+1)\bigg(\nu  + \mu (|\frac{\partial \tau}{\partial
z}|^2-|\frac{\partial \tau}{\partial \bar z} |^2) \bigg); \quad k=0,1,2, \cdots ,$$ where we have set $B:=B^{\nu,\mu}_\tau$ for simplicity.
 This is a well common fact for Schr\"odinger operators with constant magnetic field, see \cite{Fock28,FerapontovVeselov01} for example.
 Moreover, the concrete spectral analysis of $\La$, acting
on the free Hilbert space $L^2(\C; d\lambda)$, can be deduced easily by considering the transformation
\eqref{GaugeTransformation}, i.e., $
[\mathcal{W}^{\nu,\mu}_\tau  (f)](z):= e^{i \varphi^{\nu,\mu}_\tau(z)}f(z).$ It is an unitary and isometric transformation from the Hilbert
 space $L^2(\C; d\lambda)$ onto itself and intertwines the operators $\La$ and $\mathfrak{L}^{B}$ as given by
\eqref{ExplicitAutoLap} and \eqref{LandauHamiltonian}, respectively. Precisely, we have
\begin{eqnarray}\label{Intertwining} \mathcal{W}^{\nu,\mu}_\tau \La = \mathfrak{L}^{B}
\mathcal{W}^{\nu,\mu}_\tau . \end{eqnarray}
This can be checked using \eqref{Super4} and \eqref{Super5} combined with \eqref{Super1}, or also using the geometric realization of $\La$ and $\mathfrak{L}^{B}$ by writing $\nabla^{\nu,\mu}_\tau=d+i\theta^{\nu,\mu}_\tau$ as $\nabla^{\nu,\mu}_\tau= e^{-i \varphi^{\nu,\mu}_\tau}(d+i\theta^{B}) e^{i \varphi^{\nu,\mu}_\tau}$.
 Consequently, we state the following result

\begin{proposition} \label{SpecProperties} Let $A^{\nu,\mu}_{\tau;k}(\C)=\set{F\in L^2(\C; d\lambda);\quad \La F = E_k F}$ be the $L^2$-eigenspace of $\La$ associated with the eigenvalue $E_k=B(2k+1).$

  i)  We have the following orthogonal decomposition
$L^2(\C;e^{-B |z|^2} d\lambda)=\bigoplus_{k=0}^\infty A^{\nu,\mu}_{\tau;k}.$

  ii) Let $\psi^{\nu,\mu}_\tau(z,w):=\varphi^{\nu,\mu}_\tau(z)-\varphi^{\nu,\mu}_\tau(w)$. Then the eigenprojector kernel of $A^{\nu,\mu}_{\tau;k}$ satisfies the invariance property
\begin{equation*}
K^{\nu,\mu}_{\tau;k}(z,w)=
e^{-i(\psi^{\nu,\mu}_\tau(z,w)-\psi^{\nu,\mu}_\tau(g\cdot z,g\cdot w) )} e^{2i B\Im (z-w,g^{-1}.0)}K^{\nu,\mu}_{\tau;k}(g.z,g.w)
\end{equation*}
 and is given explicitly by
\[K^{\nu,\mu}_{\tau;k}(z,w)= \frac{2B}{\pi} e^{-i\psi^{\nu,\mu}_\tau(z,w)}
e^{2i B\Im \scal{z,w}} e^{-B|z-w|^2} L_k(2B  |z-w|^2)
,\] where $L_k(x)=L_k^0(x)$ denotes the usual Laguerre polynomial.
 \end{proposition}

\begin{proof} This follows by making use of \eqref{Intertwining} and Proposition 2.2 in \cite{GI-JMP08}. In fact, if $K^{\nu,\mu}_{\tau;k}(z,w)$ is a  reproducing kernel for $\La$, then $e^{-i(\varphi^{\nu,\mu}_\tau(z)-\varphi^{\nu,\mu}_\tau(w))} K^{\nu,\mu}_{\tau;k}(z,w)$ is a reproducing kernel for $\mathfrak{L}^{B}$.\end{proof}

The investigation of the spectral properties of $\La$ acting on $ \mathcal{M}^{\nu,\mu}_{\tau}(\C)$ follows now easily using Theorem \ref{Thm2} together with the fact \eqref{Intertwining}.
Indeed, we have
 $$ \mathcal{W}^{\nu,\mu}_\tau (\mathcal{E}^{\nu,\mu}_{\tau;k}) = \mathcal{E}^{B}_{k}, $$
 where for  every fixed
positive integer $k=0,1,2, \cdots$, the space $\mathcal{E}^{\nu,\mu}_{\tau;k}$ (resp. $\mathcal{E}^{B}_{k}$)
is the space of all eigenfunctions of $\La$ (resp. $\mathfrak{L}^{B}$) in $\mathcal{M}^{\nu,\mu}_{\tau}(\C)$ (resp. $\mathcal{F}^{B}_{\Gamma,\chi_\tau}(\C)$) corresponding to
the eigenvalue $E_k=B (2k+1)$, i.e.,
\begin{equation}\label{M-Eigenspace}\mathcal{E}^{\nu,\mu}_{\tau;k}:= \set{F\in
\mathcal{M}^{\nu,\mu}_{\tau}(\C); \quad
\La F = B(2k+1)F} \end{equation} and \[ \mathcal{E}^{B}_{k}:= \set{F\in
\mathcal{F}^{B}_{\Gamma,\chi_\tau}(\C); \quad
\mathfrak{L}^{B} F = B(2k+1) F}.\]
Thus,  from the dimensional formula established in \cite{GI-JMP08}, one obtains
the following

\begin{proposition} \label{DimForm}
 The  $ \mathcal{M}^{\nu,\mu}_{\tau}(\C)$-eigenspaces
$\mathcal{E}^{\nu,\mu}_{\tau;k},$ \eqref{M-Eigenspace},
 are finite dimensional spaces, whose dimension is given explicitly by
$\dim \mathcal{E}^{\nu,\mu}_{\tau;k} = (2B^{\nu,\mu}_\tau/\pi) Area(\C/\Gamma).$
\end{proposition}

For illustration, we consider the particular case of the equivariant
map $\tau_{_{\mathfrak{h}}}$, $\tau_{_{\mathfrak{h}}}(z)=\mathfrak{h} \cdot z$, derived from inner $G$-homomorphism
$\rho_{_{\mathfrak{h}}}(g) :=\mathfrak{h}g \mathfrak{h}^{-1}$;
$\mathfrak{h}=(\alpha_{_{\mathfrak{h}}},\beta_{_{\mathfrak{h}}})
\in \mathbf{T}\rtimes \C$. Thus, by solving Equation \eqref{Existence1stDiffEq},
 we get
\begin{lemma}\label{LemExpVarphi}
The function $\varphi^{\nu,\mu}_\tau$  is given explicitly by $\varphi^{\nu,\mu}_\tau(z)=\varphi^{\Nu}_\X(z):= - 2 \Im\scal{z,\X}$, where $\Nu=\nu+\mu\in \R$ and $\X=\mu\bar\alpha\beta\in\C$.
\end{lemma}
\noindent Hence, the corresponding magnetic Schr\"odinger operator acts on the Hilbert space $
L^2(\C;d\lambda)$ by
\begin{align}\label{ExplicitHammond}
 -\frac{\partial^2 }{\partial z \partial \bar z}   - \big( [\Nu z
 +\X]\frac{\partial }{\partial  z} - [\Nu\bar z  +\bar\X] \frac{\partial }{\partial \bar z} \big)
+ | \Nu z+ \X|^2  = : \mathfrak{L}_\X^\Nu,
\end{align}
 and is unitary equivalent to the twisted Laplacian $\Laa$ (\ref{LandauHamiltonian}).
Moreover, $J^{\nu,\mu}_{\rho_{\mathfrak{h}},\tau_{\mathfrak{h}}}$ reads simply
$ J^{\nu,\mu}_{\rho_{\mathfrak{h}},\tau_{\mathfrak{h}}}(\gamma;z)= e^{-2i\Nu\Im\scal{z,\gamma}}e^{-2i\Im\scal{\X,\gamma}}$. Hence, under the assumption $
 \Nu\Im \scal{\gamma,\gamma'}\in \pi \Z \label{RDC2}$, the space
of mixed $\Gamma$-automorphic forms
$\mathcal{M}^{\nu,\mu}_{\tau_{\mathfrak{h}}}(\C)$ is nontrivial and
 reduces to the space
$\mathcal{F}^{\Nu}_{\Gamma,\chi_{\mathfrak{h}}} (\C): =\set{F;  ~ \mathcal{C}^\infty,  \quad F(z+\gamma)=
\chi_{\mathfrak{h}}(\gamma)j^{-\Nu}(\gamma, z)  F(z)}$
 of Landau automorphic functions of weight $\Nu$ with
the multiplier $\chi_{\mathfrak{h}}(\gamma)=e^{-2i  \Im \scal{\X,\gamma}}$, which here is a character. Hence, one get easily that
$ \mathcal{F}^\Nu_{\Gamma;1}(\C)=\mathcal{W}_{\X}(\mathcal{F}^\Nu_{\Gamma;\chi_{\mathfrak{h}}}(\C))$ which is in accordance with Theorem \ref{Thm2}. Here $[\mathcal{W}_{\X}(f)](z):= e^{2 i\Im\scal{z,\X}}f(z)$ is the transformation in \eqref{GaugeTransformation}.
In this case, the previous results (Propositions \ref{SpecProperties} and \ref{DimForm}) read simply

\begin{corollary} \label{SpecProperties1} Let $K^{\Nu}_{\X;k}(z,w)$ denotes the eigenprojector kernel of the $L^2$-eigenspace
$A^{\Nu}_{\X;k}$ of $\mathfrak{L}_\X^\Nu$ corresponding to the eigenvalue $E_k=\Nu (2k+1)$
and $\mathcal{E}^k_{\Gamma;\Nu}(\C)$ the eigenspace of
 of $\mathfrak{L}_\X^\Nu$ in $\mathcal{M}^{\nu,\mu}_{\tau_{\mathfrak{h}}}(\C)$ associated to $E_k$. Then, we have

i)$K^{\Nu}_{\X;k}(z+b,w+b)= e^{2i \Nu\Im \scal{z-w,b}}K^{\Nu}_{\X;k}(z,w)$ for every $z,w,b\in \C$.

 ii) $K^{\Nu}_{\X;k}(z,w)= \frac{2\Nu}\pi e^{i\Im \scal{z-w,\X}}   e^{2i\Nu\Im \scal{z,w}}
e^{\Nu |z-w|^2} L_k(2\Nu |z-w|^2) $ for all $z,w\in \C$.

iii) 
$\dim
\mathcal{E}^k_{\Gamma;\Nu} = (2\Nu/\pi) Area(\C/\Gamma).$
\end{corollary}

\section{Suggested polynomials}

According to the supersymmetric relationships
\eqref{Super1}, the corresponding eigenfunctions can be generated by iterating the ground states via the creator operator $\widetilde{A^{\nu,\mu}_\tau}=  -{\partial}/{\partial z} + \overline{S^{\nu,\mu}_\tau}$. In fact, they are given by
 $ \widetilde{A^{\nu,\mu}_\tau}^m(\psi); $ $m=0,1,2, \cdots,$ where $\psi$ annihilates $A^{\nu,\mu}_\tau$. Thus, in view of Lemma \ref{Gauge}, we see that $\psi$ belongs to the space spanned by the functions $z^n e^{-B |z|^2 -i \varphi^{\nu,\mu}_\tau(z)}$; $n=0,1,2, \cdots$.
 Therefore, it follows that the free eigenspace of $\La$ corresponding to the eigenvalue $E_m=B(2m+1)$ is generated by the eigenfunctions \begin{align}\psi^{\nu,\mu;\tau}_{m,n}(z)  = \widetilde{A^{\nu,\mu}_\tau}^m(z^n e^{-B |z|^2 -i \varphi^{\nu,\mu}_\tau(z)})
  =  e^{-B |z|^2} [(\mathcal{W}^{\nu,\mu}_\tau)^{-1} \mathcal{H}^{B}_{m,n}](z),
 \end{align}
where $\mathcal{H}^{B}_{m,n}$ are the complex Hermite polynomials \cite{Shigekawa87,IntInt06} given through \begin{align*} \mathcal{H}^{B}_{m,n}(z):= e^{2B |z|^2}\frac{\partial^m}{\partial z^m} \big(z^n e^{-2B |z|^2}  \big) .\end{align*}

Now assume that $\tau$ is an equivariant map such that $R_\tau:= \tau \Delta \bar \tau - \bar \tau \Delta \tau =0$. In this case,
the corresponding $\varphi^{\nu,\mu}_\tau$ is harmonic, for $
\Delta \Im (\phi) = \Delta \varphi^{\nu,\mu}_\tau = -i\mu R_\tau  $
which follows from \eqref{SolFormDiffEq} together with the differentiation of
both sides of \eqref{1stDiffEq} with respect to $z$.
 Hence, we have $\varphi^{\nu,\mu}_\tau = 2\Im h = -i(h-\bar h)$ for some given holomorphic function $h=h^{\nu,\mu}_\tau$.
Whence, one can rewrite the annihilator and creator operators like
 $$  A^{\nu,\mu}_\tau = e^{-B |z|^2 +\bar h}
\frac{\partial}{\partial \bar z}  e^{B |z|^2- \bar h} \quad \mbox{ and } \quad
\widetilde{A^{\nu,\mu}_\tau} = e^{B |z|^2 - h}
\frac{\partial}{\partial \bar z}  e^{-B |z|^2 + h} .$$
Therefore, the functions
$$ \widetilde{\psi^{\nu,\mu;\tau}_{m,n}}(z) =  e^{-B |z|^2 + \bar h} \mathcal{G}^{B;h}_{m,n}(z)$$
generate also the free eigenspace of $\La$ with $E_m$ as associated eigenvalue, where the suggested polynomials $\mathcal{G}^{B;h}_{m,n}$ are given by
\begin{equation}
\mathcal{G}^{B;h}_{m,n}(z)= e^{2B |z|^2 - h}\frac{\partial^m}{\partial z^m} \big(z^n e^{-2B |z|^2 + h}  \big).
\end{equation}
This furnish then a class of generalized complex Hermite polynomials. Their explicit study requires explicit expression of $h$ and therefore of $\varphi^{\nu,\mu}_\tau$.

As example of equivariant pairs satisfying $R_\tau=0$, we reconsider
$(\rho_{_{\mathfrak{h}}},\tau_{_{\mathfrak{h}}})$ for which $\varphi^{\nu,\mu}_\tau=\varphi^{\Nu}_\X$ is given by Lemma \ref{LemExpVarphi}. The corresponding polynomials constitute a class of generalized complex Hermite polynomials $\mathfrak{G}^{\Nu}_{m,n}(z,\bz|\X)$ parameterized by $\X\in
\C$ and reduces to the usual complex Hermite polynomials  $\mathcal{H}^{B}_{m,n}$ when $\X=0$.  Some
related basic properties of $\mathfrak{G}^{\Nu}_{m,n}(z,\bz|\X)$ are studied and exposed in \cite{GCHP}.\\

  {\bf\it Acknowledgements.} The author acknowledges the financial support of the  Arab Regional Fellows
         Program during the academic year 2006-2007. He is indebted to Professor A. Intissar for valuable discussions and encouragements.
           A part of this work was prepared when he was a visiting member
          of CAMS, AUB.  He would like to thank all members of the center for their hospitality.


\begin{thebibliography}{9}
\bibitem{Abdulali95}  Abdulali S.,  Conjugates of
strongly equivariant maps. {\it Pacific J. Math.} 165 (1994), no. 2,
207--216.

\bibitem{Agmon95} Agmon S., Topics in spectral theory of Schrödinger operators on non-compact Riemannian manifolds with cusps.
Mathematical quantum theory. II. Schrödinger operators (Vancouver, BC, 1993), 1--33, CRM Proc.
Lecture Notes, 8, Amer. Math. Soc., Providence, RI, 1995.

\bibitem{AIM00} Askour N.; Intissar A.; Mouayn Z., Explicit formulas for reproducing kernels of generalized Bargmann spaces of $C\sp n$.
 {\it J. Math. Phys.} 41 (2000), no. 5, 3057--3067.

\bibitem{CFKS87} Cycon H. L.; Froese R. G.; Kirsch W.; Simon B. Schrödinger operators with application to quantum mechanics and global geometry.  Springer-Verlag, Berlin, 1987.

\bibitem{FerapontovVeselov01}
Ferapontov E. V.; Veselov A. P., Integrable Schrödinger operators with magnetic fields: factorization method on curved surfaces.
{\it J. Math. Phys.} 42 (2001), no. 2, 590--607.

\bibitem{Fock28}  Fock V., Bemerkung zur Quantelung des harmonischen Oszillators im Magnetfeld. Z.
Phys, 47 (1928) 446-448.

\bibitem{GCHP}  Ghanmi A., A class of generalized complex Hermite
polynomials. {\it J. Math. Anal. Appl.} 340 (2008) 1395-1406.
(arxiv:0704.3576).

\bibitem{GI05} Ghanmi A.; Intissar A., Asymptotic of complex hyperbolic geometry and $L^2$-spectral
analysis of Landau-like Hamiltonians. {\it J. Math. Phys.} 46, no.
3, (2005) 032107.

\bibitem{GI-JMP08} Ghanmi A.; Intissar A., Landau
automorphic forms of $\C^n$ of magnitude $\nu$. {\it J. Math. Phys.} 49, no. 8, (2008).
(arxiv:0705.1763).

\bibitem{HuntMeyer85}  Hunt B.; Meyer W., Mixed automorphic forms and invariants of elliptic surfaces.
Math. Ann. 271 (1985), no. 1, 53--80.

\bibitem{IntInt06}  Intissar A.; Intissar A.,  Spectral properties of the Cauchy transform on
$L^2(\C;e^{-|z|^2}d\lambda(z))$.  {\it J. Math. Anal. Appl.} 313 no
2 (2006) 400-418.

\bibitem{LeeLectNotes04}  Lee, M. H.,  Mixed automorphic forms, torus bundles, and Jacobi forms.
 Lecture Notes in Mathematics, 1845. Springer-Verlag, Berlin, 2004.

\bibitem{Matsumoto96} Matsumoto H.; Ueki N., Spectral analysis of Schrödinger operators with magnetic fields. {\it J. Funct. Anal.} 140 (1996), no. 1, 218--255.

\bibitem{SatakeBook80} Satake I., Algebraic Structures of Symmetric
Domains, Publ. Math. Soc. Japan, vol. 14, Princeton Univ. Press,
1980.

\bibitem{Shigekawa87}   Shigekawa I.,  Eigenvalue problems for the
Schr\"odinger operators with magnetic field on a compact Riemannian
manifold. {\it J. Funct. Anal.} 75  (1987) 1, 92-127.

\bibitem{Shubin01} Shubin M., Essential self-adjointness for semi-bounded magnetic Schr\"odinger operators on non-compact manifolds.
 {\it J. Funct. Anal.} 186 (2001), no. 1, 92--116.

\bibitem{Thangavelu93}  Thangavelu S., Lectures on Hermite and Laguerre expansions.
        Mathematical Notes, 42. Princeton, NJ, 1993.

\bibitem{Wong98} Wong M. W.,
  Weyl transforms. Universitext. Springer-Verlag, New York, 1998.

\end{thebibliography}
\end{document}